\documentclass[11pt]{article}
\usepackage{amssymb}
\usepackage{latexsym}
\newtheorem{thm}{Theorem}[section]
\newtheorem{prop}[thm]{Proposition} \newtheorem{lemma}[thm]{Lemma}
\newtheorem{cor}[thm]{Corollary} 
 \newtheorem{rmk}[thm]{Remark}

\newcommand {\pf}{\noindent{\bf Proof.}\ }
\newcommand{\complex}{{\mathbb C}}
\newcommand{\reals}{{\mathbb R}}

\newcommand{\cald}{{\cal D}}
\newcommand{\cale}{{\cal E}}
\newcommand{\calf}{{\cal F}}

\newcommand{\calo}{{\cal O}}

\newcommand{\calu}{{\cal U}}

\newcommand{\qed}{\begin{flushright} $\Box$\ \ \ \ \ \end{flushright}}

\newcommand{\frakg}{\mathfrak{g}}

\newcommand{\frakk}{\mathfrak{k}}

\newcommand{\fraks}{\mathfrak{s}}
\newcommand{\frakt}{\mathfrak{t}}
\newcommand{\fraksp}{\mathfrak{sp}}
\newcommand{\Ad}{\mathrm{Ad}}
\newcommand{\arrows}{\,\lower1pt\hbox{$\longrightarrow$}\hskip-.24in\raise2pt
             \hbox{$\longrightarrow$}\,}
\title{{\bf Poisson Geometry of Discrete Series Orbits, and Momentum
             Convexity for Noncompact Group Actions}}
\author{Alan
Weinstein\thanks{Research partially supported by NSF Grant DMS-99-71505.}
\\Department of Mathematics\\ University of California\\ Berkeley, CA
94720 USA\\ {\small(alanw@math.berkeley.edu)}}
\begin{document}
\setlength{\baselineskip}{15pt}

\maketitle

\section{Introduction}
\label{sec-intro}
The main result of this paper is a convexity theorem for momentum
mappings $J:M\to \frakg^*$ of certain hamiltonian actions of noncompact
semisimple Lie groups.  The image of $J$ is required to fall within a
certain open subset $\cald$ of $\frakg^*$ which corresponds roughly via
the orbit method to the discrete series of representations of the
group $G$.  In addition, $J$ is required to be proper as a map from $M$
to $\cald$.  

A related but quite different convexity theorem for noncompact groups
may be found in  \cite{ne-convexity}.

Our result is a first attempt toward placing momentum convexity
theorems in a Poisson-geometric setting.  A momentum mapping is a
Poisson mapping\footnote{All momentum mappings in this paper will be
coadjoint-equivariant, hence Poisson.}  to the dual of a Lie algebra,
but it takes more than the Poisson structure on the target manifold
even to formulate a convexity theorem.  As we explained last year's
Conference Mosh\'e Flato Proceedings \cite{we:linearization}, it seems
that a proper symplectic groupoid is the right extra structure to put
on the target.  This leads us to focus our attention on $\cald$, which
is the largest subset of $\frakg^*$ on which the (co)adjoint action is
proper, so that the subgroupoid $G\times \cald$ of the symplectic
groupoid $T^*G=G\times\frakg^*$ is a proper groupoid. It also motivates our
interest in ``relative'' convexity theorems as described in Section
\ref{sec-actions}.

$\cald$ may be defined most simply as the set of elements in $\frakg^*$
which have compact coadjoint isotropy group.   To prove that $\cald$
is open and that the action of $G$ on it is proper, we give several
other characterizations of this subset.  These involve
 the identification of
$\frakg^*$ with $\frakg$ via the Killing form and the Cartan
decomposition with respect to a maximal compact subgroup $K$ of $G$.
In particular, elements of $\cald$ correspond to Lie algebra
elements for which the corresponding vector field on the symmetric
space $G/K$ has a nondegenerate singular point.  
We also give a description of $\cald$ which shows that, as a Poisson
manifold, it is Morita equivalent, in a sense close to that of Xu 
\cite{xu:morita poisson}, to an open subset of $ \frakk^*$.
This enables us to reduce the convexity theorem to Kirwan's result
\cite{ki:convexity} for the compact case.

This paper is in some sense a companion to an earlier one on the
principal series \cite{we:principal}.  I realize with hindsight that
both papers owe a considerable amount to the paper
\cite{gu-st:frobenius} by Guillemin and Sternberg.

I would like to thank Sam Evens, Helmut Hofer, Yael Karshon, Frances Kirwan,
Toshiyuki Kobayashi, Eugene Lerman, Ken Meyer,
Wilfried Schmid, Michele Vergne, Jonathan Weitsman, and Ping Xu for helpful
comments, and Harold Rosenberg for his hospitality at the 
Institut de Math\'e\-ma\-tiques de Jussieu.  Finally, I thank
the organizers of CMF 2000 for their invitation to speak on this work,
particularly Daniel Sternheimer, whose encouragement stimulated me to
finish the proof of the convexity theorem during the long voyage from
Berkeley to Dijon.

\section{The proper part of the (co)adjoint action}
\label{sec-proper}
Let $G$ be a noncompact, semisimple reductive Lie group, $K$ a
fixed maximal compact subgroup.  Our main goal will be to relate the
coadjoint representations of $G$ and $K$. 

The Killing form of
$\frakg$ leads to the orthogonal (Cartan) decomposition $\frakg = \frakk \oplus
\fraks$, with the Killing form negative definite on $\frakk$ and
positive definite on $\fraks$.  This form and its restriction to
$\frakk$ (which is not in general the same as the Killing form of
$\frakk$) give equivariant identifications of the Lie
algebras of $G$ and $K$ with their duals.  For the rest of this
section, we will work in the Lie algebras
themselves.  In Section \ref{subsec-poisson}
we will pass to the duals and consider their
Poisson structures.

\subsection{Stable and strongly stable elements of the Lie algebra}

For $\mu \in \frakg$, its adjoint isotropy group $G_\mu$ is also the
centralizer of the closure $T_\mu$ of the 1-parameter subgroup
generated by $\mu$.  We shall say that $\mu$ is {\bf stable} if
$T_\mu$ is compact (in which case it is a torus) and {\bf strongly
stable} if $G_\mu$ is compact.  (The origin of this
terminology is explained
in Remark \ref{rmk-strongly} below.)  We denote the set of
all strongly stable elements of $\frakg$ by $\cald$.

Since every compact subgroup of $G$ is conjugate to a subgroup of $K$
(\cite{he:differential}, Thm.2.1, Ch.~VI), 
and since $\mu$ is in the Lie algebra of $T_\mu$, the adjoint
orbit of every stable $\mu$ intersects $\frakk$.  Since elements of
$\frakk$ are obviously stable, we conclude that the set of stable
elements of $\frakg$ is equal to the saturation $\Ad_G \frakk$ of
$\frakk$ by the adjoint action of $G$, and that $\cald$ is a subset
thereof.  In fact $\cald = \Ad_G \cale$, where $\cale$ is the
$K$-invariant subset $\cald \cap \frakk$ of $\frakk$.

The following proposition characterizes the elements of $\cale$ and
will lead to a description of $\cald$.

\begin{prop}
\label{prop-cale}
For $\mu \in \frakk$, the following conditions are equivalent:

(1) $\mu \in \cale$;

(2) $\frakg_\mu \cap \fraks = \{0\}$;

(3) the endomorphism $\mu_\fraks$ of $\fraks$ given 
by the infinitesimal adjoint representation is a nonsingular map;

(4) the vector field $\mu_S$ on the symmetric space $S=G/K$ given by
    the infinitesimal action of $\mu$ has a nondegenerate zero at the
    coset $eK$;

(5) $eK$ is the only zero point of $\mu_S$.
\end{prop}
\pf  (1) $\Rightarrow$ (2): The 1-parameter subgroup generated by a
nonzero element $\mu$ of $\fraks$ goes to infinity in $G$, so it cannot be
contained in a compact subgroup.  Thus $\mu$ is unstable and hence
cannot belong to $\cale$.

(2) $\Leftrightarrow$ (3): An element $w$ of $\fraks$ belongs to
    $\frakg_\mu$ if and only if $[\mu,w]=0$.  

(3) $\Leftrightarrow$ (4): The tangent space to $G/K$ at $eK$ is
    isomorphic to $\frakg/\frakk$.  The linearization of $\mu_S$ at
    $eK$ is the adjoint action of $\mu$ on $\frakg/\frakk$,
    which is equivalent to its action on $\fraks$.

(4) $\Leftrightarrow$ (5): The exponential map gives a $K$-equivariant
    diffeomorphism from $\fraks$ to $S$ (\cite{he:differential},
    Thm.~1.1, Ch.~VI), so the vector field $\mu_S$
    is equivalent to its linearization at $eK$.  

(5) $\Rightarrow$ (1) Any element of $G_\mu$ leaves the zero
set of $\mu_S$ invariant.  If this zero set reduces to the single point $eK$,
then $G_\mu$ must be contained in the isotropy group $K$ and is
therefore compact.
\qed

\begin{cor}
\label{cor-eopen}
$\cale$ is an open subset of $\frakk$.
\end{cor}
\pf
The statement follows immediately from (3) or (4) in the proposition above.
\qed

We can now characterize the strongly stable elements of $\frakg$ in
several ways.

\begin{prop}
\label{prop-cald}
For $\mu \in \frakg$, the following conditions are
equivalent.\footnote{Some numbers are omitted to make the remaining
numbers consistent with those in Proposition  \ref{prop-cale}.}

(1) $\mu \in \cald$;

(4) the vector field $\mu_S$ on the symmetric space $S=G/K$ given by
    the infinitesimal action of $\mu$ has a nondegenerate zero at some point;
    
(5) the vector field  $\mu_S$ is zero at exactly one point of $S$;

(6) $\mu$ belongs to the Lie algebra of exactly one maximal compact
    subgroup of $G$;

(7) $\mu$ is an interior point of the set $\Ad_G \frakk$ of stable elements.
\end{prop}
\pf
(1) $\Leftrightarrow$ (4): If $\mu\in \cald$, then $\Ad_g \mu \in
\cale$ for some $g 
\in G$.  The action of $g$ on $S$ gives an equivalence between the
vector fields $(\Ad_g \mu)_S$ and $\mu_S$.  By (4) of Proposition
\ref{prop-cald}, $(\Ad_g \mu)_S$ has a nondegenerate zero at a point
of $S$, hence so does $\mu_S$.  Conversely, if $\mu_S$ has a zero at
some point, then $\Ad_g \mu \in
\cale$ for some $g\in G$.  If the zero of $(\Ad_g \mu)_S$ is
nondegenerate, so is that of $\mu_S$, in which case $Ad_g \mu \in \cale$,
and hence $\mu \in \cald$.  

(1) $\Leftrightarrow$ (5): The argument is just like the one which
    showed that (1) $\Leftrightarrow$ (4).  Or one can show that (4)
    $\Leftrightarrow$ (5) by using the riemannian exponential map at
    the zero point.

(5) $\Leftrightarrow$ (6): The maximal compact subgroups are the
    isotropy groups of the points of $S$.

(1) $\Leftrightarrow$ (7): By Corollary \ref{cor-dopen} below, $\cald$
    is open in $\frakg$.  Since $\cald \subset \Ad_G \frakk$, all
    points of $\cald$ are interior points of $\Ad_G \frakk$.
    Conversely, suppose that $\mu$ is an interior point of $\Ad_G
    \frakk$.  Then we have some $\nu=\Ad_g \mu \in \frakk$, and we
    need to show that $\nu$, which is also an interior point of $\Ad_G
    \frakk$, belongs to $\cale$.  If it did not, there would be a
    nonzero $w\in \fraks$ with $[\nu,w]=0$.  But then, for all real $t$,
    $\nu+tw$ would be unstable, and hence $\nu$ could not be an
    interior point of $\Ad_G
    \frakk$.
\qed

\begin{cor}
\label{cor-dopen}
$\cald$ is an open subset of $\frakg$.
\end{cor}
\pf
The statement follows directly from (4) in the proposition above,
since the property of having a nondegenerate zero point is stable
under small perturbations.
\qed

\begin{rmk}
\label{rmk-strongly}
{\em
Part (7) of Proposition \ref{prop-cald}
says that a stable element of $\frakg$ is
strongly stable  if and only it remains stable under small
perturbations.  This property justifies the term ``strongly stable''--
the term
 was originally used in the case where $G$ is a symplectic group
and $K$ a unitary group, elements of $\frakg$ then being linear
autonomous hamiltonian dynamical systems.  See
Section \ref{subsec-example}
below for further discussion of this example, with some references.  
}
\end{rmk}

We may now identify those groups for which $\cald$ is nonempty.  Since
every strongly stable element is stable, it is semisimple, and since
$\cald$ is open, it contains regular semisimple elements.  The
stabilizer of such an element is a compact Cartan subgroup of $G$
(\cite{he:differential} Thm.~3.1, Ch.~III).

Conversely, if $G$ has a compact Cartan subgroup $T$, the stabilizer
of any regular element of $\frakt$ is equal to $T$, hence compact, so
this element belongs to $\cald$.  Thus we have shown:

\begin{prop}
\label{prop-nonempty}
The set $\cald \subset \frakg$ of strongly stable elements is nonempty
if and only if $G$ contains a compact Cartan subgroup, i.e.\ if and only if
$G$ has the same rank as its maximal compact subgroup $K$.  
\end{prop}

\begin{rmk}
\label{rmk-discrete}
{\em We note that the condition $\mathrm{rank}~K=\mathrm{rank}~G$ is
precisely the condition for $G$ to admit discrete series
representations \cite{ha:discrete2}
i.e.\ irreducible subrepresentations of its
left regular representation.  In fact, the discrete series
representations correspond via character theory or
the orbit method to certain orbits in
$\cald$. (See Remarks \ref{rmk-discrete2} and \ref{rmk-dual} below.)
}
\end{rmk}

\subsection{Characterization in terms of roots}

To give a more concrete description of the elements of $\cale$, and
hence those of $\cald$, we choose a maximal torus $T\subset K$
(which is also a Cartan subgroup of $G$).  Since $\Ad_K\frakt$ is all
of $\frakk$, we have $\cale = \Ad_K \calf$, where $\calf = \cale \cap
\frakt$.  Then $\cald = \Ad_G \calf$.

To identify the elements of $\calf$, we recall that the
adjoint action of $T$ on $\frakk$ leaves invariant the splitting
$\frakk\oplus\fraks$, so that each root in $\frakt^*$ can be designated
as either ``compact'' or ``noncompact'' according to whether the
corresponding eigenvector lies in $\frakk$ or $\fraks$.  Then the
following characterization of $\calf$ follows immediately from (3) of
Proposition \ref{prop-cale}.

\begin{prop}
\label{prop-calf}
$\cald = \Ad_G \calf$, where
$\calf \subset \frakt$ is the complement of the zero-hyperplanes of
the noncompact roots.  In particular, $\calf$ is dense in $\frakt$ and
is the disjoint union of finitely many convex open subsets.
\end{prop}

\begin{rmk}
\label{rmk-discrete2}
{\em
According to Harish-Chandra \cite{ha:discrete2}, 
the discrete series representations are parametrized 
by the subset of
$\cald$
consisting of those (integral) orbits
whose intersection with $\frakt$
lies in the complement of the zero-hyperplanes of {\em
all} the roots.  The remaining integral orbits in $\cald$ correspond
to ``limits of discrete series representations.''
}
\end{rmk}

The density of $\calf$ in
$\frakt$ implies that the stable and strongly stable 
elements of $\frakg$ have the same closure, i.e.:

\begin{cor}
\label{cor:closures}
$\overline{\cald}=\overline{\Ad_G \frakk}$.
\end{cor}

The example of $SL(2,\reals)$ already shows that $\Ad_G \frakk$ itself
is not closed--its closure includes the nilpotent cone.

\subsection{Example: the symplectic group}
\label{subsec-example}
Among the groups having the same rank as their maximal compact
subgroups are the indefinite unitary groups $U(m,n)$, the
indefinite orthogonal groups $SO(m,n)$ for which $m$ and $n$ are not
both odd, and the (real) symplectic groups $Sp(2n)$.  In this section,
we will concentrate on the last example and will study the strongly
stable elements of $\fraksp(2n)$.

We use canonical coordinates $(q_j,p_j)$ for the symplectic structure
$\omega=\sum dq_j\wedge dp_j$ on $\reals^{2n}$.  We
 also identify $\reals^{2n}$ with $\complex^n$ by using the complex coordinates
 $z_j=q_j+ip_j$.  The maximal compact subgroup of $G=Sp(2n)$ is the
 unitary group $K=U(n)$; the Cartan subgroup $T=T^n$ consists of the
 diagonal unitary matrices.  The symmetric space $S=Sp(2n)/U(n)$ may
 be identified with the set of positive polarizations on
$\reals^{2n}$, i.e.\ the almost complex structures $J$ on
 $\reals^{2n}$ for which the bilinear form $\omega(x,Jy)$ is
symmetric and positive-definite.  

We identify elements of $\fraksp(2n)$ with linear hamiltonian vector
fields and, in turn, with the quadratic hamiltonian functions which
generate them.

The equivalence of (5) and (7) in Proposition \ref{prop-cald} tells us
that a linear hamiltonian system is stable and remains so under small
perturbations if and only if it is ``uniquely unitarizable,'' i.e.\ if
an only if it leaves invariant a unique compatible complex structure.
This type of result plays a basic role in 
I.~Segal's approach to quantum field theory.  For instance, in
\cite{se:quantization}, Segal proves in the infinite dimensional case
that a stable linear symplectic map $T$ is ``uniquely unitarizable''
if and only if $T$ and $T^{-1}$ have disjoint spectra.  

We can also use Proposition \ref{prop-calf} to 
characterize strongly stable quadratic hamiltonians.  The compact
Cartan subalgebra consists of hamiltonians of the form 
\begin{equation}
\sum_{j=1}^n \frac{\lambda_j}{2}(q_j^2+p_j^2),
\end{equation}
where the $\lambda_j$ are arbitrary real numbers whose absolute values
are the frequencies of the normal modes of oscillation for the linear
hamiltonian system.

The compact roots are the differences $\lambda_j - \lambda_k$ for all
pairs $j < k$, while the noncompact roots are the sums $\lambda_j
+\lambda_k$ for $j\leq k$.   The strong stability condition --
nonvanishing of the noncompact roots -- means that all the normal mode
frequencies are nonzero, and that, whenever there is a  simple resonance
$|\lambda_j|=|\lambda_k|$, the signs of $\lambda_j$ and $\lambda_k$ are
the same.  This criterion for strong stability is well known in the
theory of hamiltonian systems. \footnote{Strong stability, also known
as {\bf parametric stability}, was characterized in a similar way
for elements of the symplectic group, rather
than of its Lie algebra, by Krein \cite{kr:generalization} (who only announced
results), and Gelfand and Lidskii \cite{ge-li:structure}.
The same characterization (definiteness of a quadratic form on
eigenspaces) was given by Moser \cite{mo:new}.}

Note that all positive-definite and negative-definite hamiltonians are
strongly stable--these form two of the connected components of
$\cald$.

\subsection{Bundle structure; properness of the action}
We still owe the reader a proof of the property which originally
motivated our definition of $\cald$, namely the properness of
the adjoint action.  The proof will be based on the following
description of $\cald$.  

\begin{prop}
\label{prop-bundle}
$\cald$ is $G$-equivariantly isomorphic to the associated bundle
$(G\times \cale)/K$, which is an open subbundle of the homogeneous
vector bundle $(G\times \frakk)/K$.
\end{prop}
\pf
By (4) and (5) of Proposition \ref{prop-cald}, for each $\mu \in
\cald$ the vector field
$\mu_S$ has a unique zero $\phi(\mu)$ in $S=G/K$, and this zero is nondegenerate.
It follows from the implicit function theorem that $\phi$ is a smooth
mapping from $\cald$ to $G/K$.  Since the map $\mu\mapsto\mu_S$ is
$G$-equivariant, so is $\phi$.  The fibre of $\phi$ over $eK$ consists
of those $\mu\in \cald$ which generate 1-parameter subgroups fixing
$eK$, i.e.\ $\phi^{-1}(eK)=\cald\cap\frakk=\cale$.  The statement of
the proposition follows.  Concretely, we map $G \times \cale$ to
$\cald$ by $(g,\nu)\mapsto \Ad_g \nu$ and observe that the fibres of
this map are the $K$-orbits. 
\qed

We will use the following simple lemma about proper actions.

\begin{lemma}
\label{lem-proper}
If $\phi:X\to Y$ is a continuous equivariant map of $G$-ma\-ni\-folds,
and if the
action of $G$ on $Y$ is proper, so is the action on $X$. 
\end{lemma}
\pf
By the definition of properness, the action map $\alpha_Y:(g,y)\mapsto(gy,y)$
from $G\times Y$ to $Y\times Y$ is proper.  To check that the
corresponding map $\alpha_X$ for the action on $X$ is proper, we let
$K$ be an arbitrary compact subset of $X\times X$.  Equivariance
implies that $(\mathrm{Id} \times
\phi)(\alpha_X^{-1}(K))\subseteq
\alpha_Y^{-1}((\phi\times\phi)(K)).$ Since $K$ is compact
and $\phi$ is continuous, $(\phi\times\phi)(K)$ is compact; 
by the properness of $\alpha_Y$, $\alpha_Y^{-1}((\phi\times\phi)(K))$
is compact as well, hence so is its closed subset 
$(\mathrm{Id} \times
\phi)\alpha_X^{-1}(K)$.  Applying the projection to $G$, which 
is unaffected by $(\mathrm{Id} \times \phi)^{-1},$ we find that there
is a compact subset $A\subseteq G$ such that $A\times X$ contains
$\alpha_X^{-1}(K)$.  On the other hand, applying the second projection
from $K\subseteq X\times X$ into $X$, which is unaffected by
$\alpha_X$, we find that $\alpha_X^{-1}(K)$ is contained in $G\times
B$, where $B$ is a compact subset of $X$.  Thus the closed set
$\alpha_X^{-1}(K)$ is contained in $A\times B$, and so it is compact.
\qed

\begin{cor}
\label{cor-proper}
The adjoint action of $G$ on $\cald$ is proper.  
\end{cor}
\pf  The result follows from Lemma \ref{lem-proper} and
Proposition \ref{prop-bundle} once we
observe that the action map $G\times G/K \to G/K \times G/K$ is a
locally trivial fibration with compact fibres (essentially
the isotropy groups of the action).
\qed

\subsection{Poisson geometry of $\cald$}
\label{subsec-poisson}
We will now relate the Poisson geometry of $\cald \subset
\frakg^*$ to that of $\cale \subset \frakk^*$ by constructing a Morita
equivalence, or dual pair, relating them.  The
construction is based on the slice method of \cite{gu-st:convexity}, 
where an analogous
relation is established between an open subset of $\frakk*$ (the
regular elements) and the interior of the
positive Weyl chamber in $\frakt^*$.

The Killing form of $\frakg$ gives equivariant identifications of
$\frakg^*$ with $\frakg$ and of $\frakk^*$ with $\frakk$ and 
the annihilator of $\fraks$ in $\frakg^*$.  This allows us to consider
$\cale$ and $\cald$ as open Poisson submanifolds of $\frakg^*$ and
$\frakk^*$ respectively.

We will call a submanifold $N$
of a Poisson manifold $P$  {\bf cosymplectic} if the
restriction of the Poisson tensor to each conormal space of $P$ is
nondegenerate.  Equivalently, if $P$ is defined locally by the
vanishing of independent functions $f_1,\ldots f_k$, $P$ is
cosymplectic if the matrix of Poisson brackets $a_{ij}=\{f_i,f_j\}$ is
nondegenerate along $P$.  (In the language of Dirac
\cite{di:lectures}, the $f_i$ are {\em second-class constraints}.)
Another characterization is that $N$ is cosymplectic if it is
transversal to each symplectic leaf of $P$, and if its
intersection with each such leaf
is a symplectic
submanifold of the leaf.  In particular, the cosymplectic submanifolds
of a symplectic manifold are just the symplectic submanifolds.

\begin{lemma}
\label{lem-cosymplectic}
Let $N$ be a cosymplectic submanifold of the Poisson manifold $P$, and
let $\phi:Q\to P$ be a Poisson map. Then $\phi$ is transverse to $N$, and
$\phi^{-1}(N)$ is a cosymplectic submanifold of $Q$.  
\end{lemma}
\pf
Let $\phi(x)\in N$ for some $x\in Q$, and let $(f_1,\ldots, f_k)$ be
independent functions defining $N$ near $\phi(x)$.  Then the functions
$(\phi^* f_1,\ldots, \phi^* f_k)$ define $\phi^{-1}(N)$ near $x$.  The
matrix of Poisson brackets $$\{\phi^* f_i,\phi^* f_j\}(x) =
\phi^*\{f_i,f_j\}(x)=\{f_i,f_j\}(\phi(x))$$
is nondegenerate.  It follows, first of all, that the 
functions
$(\phi^* f_1,\ldots, \phi^* f_k)$ are independent near $x$, so that
$\phi$ is transverse to $N$.  It then follows that $\phi^{-1}(N)$ is
cosymplectic.  
\qed

\begin{prop}
\label{prop-cosymplectic}
The subset $\cale\subset \frakk^*$ of strongly stable points is the set
of points at which $\frakk^*$ is cosymplectic in $\frakg^*$.  Hence
$\cale$ is a cosymplectic submanifold of $\cald$.  
\end{prop}
\pf
The conormal space to $\frakk^* \subset \frakg*$ at $\mu$ may be
identified with $\fraks \subset \frakg$, on which the Poisson tensor
is the bilinear form $(v,w)\mapsto \langle \mu,[v,w]\rangle$, where
$\langle~,~\rangle$ is the Killing form.  By invariance of the Killing
form, this can also be written as $<[\mu,v],w>$, so the Poisson tensor is
nondegenerate exactly when the action of $\mu$ on $\fraks$ is
nonsingular.  Now apply (3) of Proposition \ref{prop-cale}.
\qed

We shall identify the symplectic manifold $T^*G$ with $G\times
\frakg^*$ by right translations.  Then we have:
\begin{cor}
\label{cor-cosymplectic}
The product $M=G\times \cale \subset G\times \frakg$ is a symplectic
submanifold of $T^*G$.

The action of $G$ on $M$ by right translations
($g\cdot h = hg^{-1}$) is hamiltonian
with momentum map equal to $\psi:(g,\mu)\mapsto -\Ad _g \mu$, which
is a Poisson submersion from $M$ to $\cald$.  $M$ is also
invariant under the hamiltonian action of $K$ by left translations,
for which the momentum map is right translation, i.e.\ the projection
from $M=G\times \cale$ to $\cale$.  The two momentum maps form a
symplectic dual pair relating the Poisson manifolds $\cale$ and $\cald$.  
\end{cor}
\pf  The momentum map for the cotangent lift of the action of $G$ on
itself by right translations is the negative of the left translation
map.  In our right trivialization, this appears as $\psi:(g,\mu)\mapsto
-\Ad _g \mu$. It follows from Lemma \ref{lem-cosymplectic}
and
Proposition \ref{prop-cosymplectic} that $M$ is
a symplectic submanifold of $T^*G$.  In general, the momentum map for
the restriction of a hamiltonian action to an invariant symplectic
submanifold is the restriction of the original momentum map, since the
momentum map for the right action of $G$ on $M$ is the negative of the
projection, hence $\psi:(g,\mu)\mapsto
-\Ad _g \mu$ is a Poisson map from $M$ to $\frakg$.   It is a
submersion because the action of $G$ on $M$ is free, and its image 
is $\Ad_G \cale = \cald$.  

Invariance of $M$ under left translations by $K$ follows from the
$\Ad_K$-invariance of $\cale$.  The momentum map for the left action
of $K$ on $T^*G = G\times \frakg^*$ is given by right translation,
which is just the projection.  Again, it is a submersion because the
action of $K$ is free.  

To see that the pair of maps $\frakg\leftarrow M \rightarrow \frakk$
form a symplectic dual pair, i.e.\ that they are submersions whose
fibres have symplectically orthogonal tangent spaces, we simply note that
they are the momentum maps of commuting free actions of groups whose
dimensions add up to that of $M$.  
\qed

\begin{rmk}
\label{rmk-dual}
{\em
Assuming that
that $G$ (and hence $K$) is connected, the dual pair of Proposition
\ref{prop-cosymplectic} has
connected fibres.  By Proposition 9.2 of \cite{ca-we:geometric}, this 
gives a bijection between the symplectic leaves of $\cald$
and $\cale$.  Here, this bijection is simply given by intersecting the
leaves with $\frakk$.

An important property of dual pairs of group actions is that the orbit
space of one action is Poisson-isomorphic to the image of the momentum
map of the other.  Using this fact one way gives the uninteresting
isomorphism between $(G\times \cale)/G$ and $\cale$.  But in the other
direction we obtain an isomorphism between $(G\times \cale/K)$ and
$\cald$.  This is a symplectic version of the  description of $\cald$
as an $\cale$ bundle over $G/K$ given in \ref{prop-bundle}.  

A symplectic dual pair, i.e.\ a pair of Poisson submersions with
connected, symplectically orthogonal fibres, is called a {\bf Morita
equivalence} if the fibres are
simply-connected and if the submersions satisfy a completeness
condition.  Xu \cite{xu:morita poisson}, has shown that Morita
equivalence of Poisson manifolds implies {\bf representation
equivalence}, i.e.\ an equivalence between their symplectic
realizations.  Although in some examples the dual pair of P
\ref{prop-cosymplectic} does not have simply-connected fibres, one
may think of the representation equivalence between $\cald$ and
$\cale$ as a classical analogue of the equivalence between the
discrete series representations of $G$ and the representations of the
maximal compact subgroup $K$.
}
\end{rmk}

\section{Actions and convexity}
\label{sec-actions}

We will begin this section by defining properness of momentum maps
relative to open subsets of the dual of a Lie algebra.  After
recalling the known convexity theorems for compact group actions, we
will prove the central theorem of this paper, a convexity theorem for
certain actions of noncompact groups.

\subsection{Proper hamiltonian spaces}
\label{subsec-proper}

Let $G$ be a Lie group, and let $\calu\subseteq \frakg^*$  be a
coadjoint-invariant  open subset.  We define a {\bf hamiltonian}
$(G,\calu)$-{\bf  space} $(M,J)$ to be a symplectic manifold $M$ with
a symplectic $G$-action and a coadjoint-equivariant momentum map
$J:M\rightarrow \calu \subseteq \frakg^*.$  We shall consider $J$ as a
map to $\calu$ rather than to $\frakg^*$ and will call the
$(G,\calu)$ space {\bf proper} if $J$ is a proper mapping and if the
action of $G$ on $M$ is proper.  By Lemma \ref{lem-proper}, the second
condition follows from the first if the coadjoint action of $G$ on
$\calu$ is proper, e.g. when $G$ is semisimple and $\calu$ consists of
strongly stable elements.

\subsection{The compact case}
\label{subsec-compact}

The first convexity theorem applies to torus actions.

\begin{thm}
\label{thm-torus}
Let $\calu$ be a disjoint union of convex open subsets
 of the dual of the Lie algebra of a
torus $T$, and let $(M,J)$ be a connected,
proper, hamiltonian $(T,\calu)$-space.
Then $J(M)$ is a closed, convex locally polyhedral subset of $\calu$, and the
inverse image $J^{-1}(\mu)$ is connected for each $\mu\in\calu$.  
\end{thm}

Theorem \ref{thm-torus} was originally proved by Guillemin and
Sternberg \cite{gu-st:convexity} and Atiyah \cite{at:convexity}  in
the case where $M$ is compact and $\calu$ is all of $\frakt^*.$
Extensions to noncompact $M$ were given by Prato \cite{pr:convexity}
and, using methods of Condevaux, Dazord, and Molino
\cite{co-da-mo:geometrie},
 by Hilgert, Neeb, and Plank
\cite{hi-ne-pl:convexitycoadjoint}.

An extension of the convexity theorem to actions of nonabelian compact
groups on compact manifolds was conjectured and partially proved by
Guillemin and Sternberg \cite{gu-st:convexity2}.  The proof was
completed by Kirwan \cite{ki:convexity}.  The result was
extended to proper actions in \cite{hi-ne-pl:convexitycoadjoint}, and
the ``relative'' version, where the momentum map is proper as a map
into an open subset containing its image, was obtained by Lerman,
Meinrenken, Tolman, and Woodward \cite{le-me-to-wo:nonabelian}.

In the statement of the following theorem,
we have identified the Lie algebras with their duals by using a
bi-invariant metric on the group.

\begin{thm}
\label{thm-compact}
Let $K$ be a compact Lie group, let $\frakt_+^*$ be a positive Weyl
chamber in $\frakt^*\subseteq \frakg^*$, and let $\calu \subseteq\frakg$ be a
coadjoint-invariant open subset such that each component of 
$\calu\cap\frakt_+^*$ is convex. If $(M,J)$ is a connected, proper,
hamiltonian $(G,\calu)$-space, then $J(M)\cap \frakt_+^*$ is a closed,
convex, locally
polyhedral subset of $\calu\cap\frakt_+^*,$ and $J^{-1}(\mu)$
is connected for each $\mu\in\calu.$
\end{thm}

\subsection{The noncompact case}

We will now state and prove the main result of this paper.

\begin{thm}
\label{thm-main}
Let $G$ be a semisimple Lie group, let $\frakt_+^*$ be a positive Weyl
chamber for a maximal compact subgroup $K$ of $G$, and 
let $\calu$ be a coadjoint-invariant open subset of the set
$\cald\subset\frakg*$ of strongly stable elements such that
$\calu\cap\frakt_+^*$ is convex.
 If $(M,J)$ is a connected, proper, hamiltonian $(G,\calu)$-space, then 
$J(M)\cap \frakt_+^*$ is a closed, convex, locally polyhedral subset of
$\frakt_+^*\cap \calu,$ and $J^{-1}(\mu)$ is connected for each
$\mu\in\calu.$
\end{thm}

\pf
Just as the convexity theorem for nonabelian compact $K$ is proved by
reduction to the abelian case, our theorem for noncompact $G$ will be
proved by reduction to the compact case.  To this end, we let
$N=J^{-1}(\calu\cap\cale)  = J^{-1}(\calu\cap\frakk^*)
\subseteq M$.  Since $\cale$ is cosymplectic in $\cald$, $\calu \cap
\cale$ is cosymplectic in $\calu$; by Lemma \ref{lem-cosymplectic}, 
$J$ is transverse to
$\calu\cap\cale,$  and $N$ is a 
cosymplectic submanifold of $M$.  Thus, $N$ is a symplectic manifold,
and, since $J$ is equivariant, $N$ is $K$-invariant.   By basic facts
about the behavior
of momentum maps under restriction (to subgroups and
submanifolds), 
$J|_N:N\rightarrow \calu\cap\cale$ is an 
equivariant momentum map for the action of $K$ on $N$, making $N$ into
a hamiltonian $(K,\calu\cap\cale)$-space.  The properness of $J$
implies immediately that $J|_N$ is proper as well.

To apply Theorem \ref{thm-compact} to this $K$-space, we just need to
know that $N$ is connected.  To prove this, we compose $J$ with the
equivariant map $\phi:\cald\to G/K$ of Proposition
\ref{prop-bundle}.  Since $\cale=\phi^{-1}(eK)$, $\phi\circ J$ makes
$M$ into a bundle over $G/K$ with typical fibre $J^{-1}(\cale)=N$.
Now $G/K$ is contractible, hence simply connected, and $M$ is
connected, so $N$ must be connected as well.  It follows from Theorem
\ref{thm-compact} that $J|_N(N)\cap \frakt_+^*$ is a closed, convex,
locally polyhedral subset of $ \frakt_+^*\cap (\calu\cap\cale)=
\frakt_+^*\cap\cale,$ and that $J|_N^{-1}(\mu)$ is connected for each
$\mu\in \calu\cap\cale.$  But 
$J|_N(N)\cap \frakt_+^*$ is equal to $J(M)\cap\frakt_+^*:$ also, for
any $\nu$ in the invariant subset $\calu\subseteq \cald \subset \Ad
_G\cale,$ so $\nu=g\mu$ for some $\mu$ in $\calu\cap\cale,$ so
$J^{-1}(\nu)=gJ^{-1}(\mu)$ is connected.  This completes the proof of
the theorem.
\qed

\section{Products of coadjoint orbits and sums of positive definite
hamiltonians} 
\label{sec-products}
Let $\calo_1$ and $\calo_2$ be coadjoint orbits in $\frakg^*$.  Then
$G$ acts diagonally on the product $\calo_1 \times \calo_2$ with
momentum map given by the addition map $A:\calo_1 \times \calo_2\to
\frakg^*$.  There are certain cases in which the image of $A$ lies in
$\cald$. For instance, in the case of the symplectic groups (see
Section \ref{subsec-example}), we may let $\calo_1$ and $\calo_2$ be
orbits consisting of positive-definite hamiltonians, since the sum of
two positive functions is positive.  In this case, our convexity
theorem has the following interesting consequence in hamiltonian
dynamics.

\begin{thm}
\label{thm-sums}
For any positive-definite quadratic hamiltonian function $H$ on
$\reals^{2n}$, let $F(H)$ be the $n$-tuple $(\lambda_1,\ldots
\lambda_n)$, where $\lambda_1\leq \cdots\leq\lambda_n$ are the frequencies
of the normal modes of oscillation for the linear hamiltonian system
generated by $H$; i.e.\ $F(H)$ are the coefficients of the normal form 
$\sum_{j=1}^n \frac{\lambda_j}{2}(q_j^2+p_j^2)$ for $H$ in suitably
chosen canonical coordinates.  If $\lambda$ and $\mu$ are 
nondecreasing $n$-tuples of positive real numbers, then 
\begin{equation}
\label{eq-sums}
\{F(H_1+H_2)|F(H_1)=\lambda {\rm ~and~} F(H_2)=\mu\}
\end{equation}
is a closed, convex, locally polyhedral set.
\end{thm} 
\pf
We follow a standard approach, used already to study the spectrum of sums of
hermitian matrices with prescribed eigenvalues.  Let $\calo_{\lambda}$
and $\calo_{\mu}$ be the (co)adjoint orbits of
$\sum_{j=1}^n \frac{\lambda_j}{2}(q_j^2+p_j^2)$
and $\sum_{j=1}^n \frac{\mu_j}{2}(q_j^2+p_j^2)$
respectively.  Theorem \ref{thm-sums} will follow from Theorem
\ref{thm-main} once we have proven that the addition
operation $(H_1,H_2)\mapsto H_1 \times
H_2$ is a proper map from $\calo_{\lambda}\times\calo{\mu}$ to $\cald$.

Suppose, then, that $\{H_{1,i}\}$ and $\{H_{2,i}\}$ are sequences in
$\calo_{\lambda}$ and $\calo_{\mu}$ respectively such that
$H_{1.i}+H_{2,i}$ converges to an element $H$ of $\cald$.  We must
show that $\{(H_{1,i},H_{2,i})\}$ has a convergent subsequence.  

We will use the ``diagonal'' projection $\pi$ from $\fraksp(2n)$ to its
Cartan subalgebra $\reals^n$, which selects the coefficients of the
terms $\frac{1}{2}(q_j^2+p_j^2)$ in a quadratic hamiltonian, ignoring
the coefficients of $\frac{1}{2}(q_j^2-p_j^2)$ and all ``cross
terms''.  The map $\pi$ is the momentum map for the action
of the Cartan subgroup and, according to Proposition 2.2
in \cite{pr:convexity},
its restriction to any coadjoint orbit of positive-definite
hamiltonians is proper.   

By applying a preliminary transformation in $Sp(2n)$, we may assume
that $H$ is in the normal form $\sum_{j=1}^n
\frac{\nu_j}{2}(q_j^2+p_j^2)$.  Since $\pi$ is linear and takes
positive definite matrices to positive $n$-tuples, the sequences
$\{\pi(H_{1,i})\}$ and $\{\pi(H_{2,i})\}$ are bounded and, hence,
after we pass to subsequences, may be assumed convergent.  But now the
properness of $\pi$ on the coadjoint orbits implies that
$\{(H_{1,i},H_{2,i})\}$ has a convergent subsequence, and the theorem
is proven.
\qed
 
We close this section on a speculative note.  Is there a version of
this convexity theorem for nonlinear symplectic transformations?  To
formulate the question more precisely, we note that the fundamental
frequencies $\lambda_i$ of a quadratic hamiltonian can be interpreted
as the values of the ratio $$\rho = 2\pi\frac{\mathrm{energy}}{\mathrm
{action}}$$ for the simple periodic trajectories of the
hamiltonian dynamical system.  We may ask, then, about the possible
values of $\rho$ for periodic orbits for  
$\sum_{j=1}^n
\frac{\lambda_j}{2}(q_j^2+p_j^2)\circ \phi_1 + \sum_{j=1}^n
\frac{\mu_j}{2}(q_j^2+p_j^2)\circ\phi_2$, where $\phi_1$ and $\phi_2$
are homogeneous, but not necessarily linear, symplectic
diffeomorphisms of $\reals^{2n}\backslash \{0\}$.


\begin{thebibliography}{99}

\bibitem{at:convexity}
Atiyah, M.F., Convexity and commuting Hamiltonians, {\em Bull. London
Math. Soc.} {\bf 14} (1982), 1-15.

\bibitem{ca-we:geometric}
Cannas da Silva, A., and Weinstein, A., {\em 
Geometric Models for Noncommutative Algebras}, Berkeley
Math. Lecture Notes, Amer. Math. Soc., Providence, 1999.

\bibitem{co-da-mo:geometrie}
Condevaux, M., Dazord, P., and Molino, P., G\'eom\'etrie du moment,
Travaux du S\'eminaire Sud-Rhodanien de G\'eom\'etrie, I, {\em
Publ. D\'ep. Math.  Nouvelle S\'er. B} {\bf 88-1},
Univ. Claude-Bernard, Lyon, 1988, 131--160.

\bibitem{di:lectures}
Dirac, P. A. M., {\em Lectures on quantum mechanics}, 
Belfer Graduate School of Science, Yeshiva University, New York, 1964.

\bibitem{ge-li:structure}
Gel'fand, I. M., and  Lidskii, V. B. ,
On the structure of the regions of stability of linear canonical
systems of differential equations with periodic coefficients,
(Russian) {\em Uspehi Mat. Nauk (N.S.)} {\bf 10} (1955), 3--40. 

\bibitem{gu-st:convexity} 
Guillemin, V., and Sternberg, S., Convexity
properties of the moment mapping, {\em Invent. Math.} {\bf 67} (1982),
491-513.

\bibitem{gu-st:frobenius}
Guillemin, V., and Sternberg, S., The Frobenius reciprocity theorem
from a symplectic point of view, {\em Lecture Notes in Math} {\bf
1037} (1983), 242-256.

\bibitem{gu-st:convexity2} 
Guillemin, V., and Sternberg, S., Convexity
properties of the moment mapping II, {\em Invent. Math.} {\bf 77} (1984),
533-546.

\bibitem{ha:discrete2}
Harish-Chandra, Discrete series for semisimple Lie
 groups.~II. Explicit determination of
the characters, {\em Acta Math.} {\bf 116} (1966), 1--111.

\bibitem{he:differential}
Helgason, S., {\em Differential geometry and symmetric spaces},
Academic Press, New York (1962).

\bibitem{hi-ne-pl:convexitycoadjoint}
Hilgert, J., Neeb, K.-H., and Plank, W., Symplectic
convexity theorems and coadjoint orbits, {\em Compositio Math.} {\bf
94} (1994), 129--180. 

\bibitem{ki:convexity} 
Kirwan, F., Convexity
properties of the moment mapping III, {\em Invent. Math.} {\bf 77} (1984),
547-552.

\bibitem{kr:generalization}
Krein, M. G., A generalization of some investigations of
A. M. Lyapunov on linear differential equations with periodic
coefficients (Russian) {\em Doklady Akad. 
Nauk SSSR (N.S.)} {\bf 73}, (1950), 445--448.  

\bibitem{le-me-to-wo:nonabelian}
Lerman, E., Meinrenken, E., Tolman, S., and  Woodward, C.,
Nonabelian convexity by symplectic cuts, {\em  Topology} {\bf 37}
(1998), 245--259.

\bibitem{mo:new}
Moser, J.,
New aspects in the theory of stability of Hamiltonian systems,
{\em Comm. Pure Appl. Math.} {\bf 11} (1958), 81--114. 

\bibitem{ne-convexity}
Neeb, K.-H.,
Convexity properties of the coadjoint action of non-compact Lie
groups, {\em Math. Ann.} {\bf 309} (1997), 625--661. 

\bibitem{pr:convexity}
Prato, E., Convexity properties of the moment map 
for certain non-compact manifolds, {\em Comm. Anal. Geom.} {\bf 2} (1994),
267--278. 

\bibitem{se:quantization}
Segal, I., Quantization of symplectic transformations,
Mathematical analysis and applications, Part B, 
{\em Adv. in Math. Suppl. Stud.} {\bf 7b}, 
Academic Press, New York-London (1981), 749--758. 

\bibitem{we:principal}
Weinstein, A., Poisson geometry of the principal series and
nonlinearizable structures, {\em J. Diff. Geom.} {\bf 25} (1987), 55--73.

\bibitem{we:linearization}
Weinstein, A., Linearization problems for Lie algebroids and Lie
groupoids, {\em Lett. Math. Phys.} {\bf 52} (2000), 93-102.


\bibitem{xu:morita poisson}
Xu, P. Morita equivalence of Poisson manifolds, 
{\em Comm. Math. Phys.} {\bf 142} (1991), 493-509.
\end{thebibliography}
\end{document}